\documentclass[letterpaper]{IEEEtran} 

\usepackage{cite}
\usepackage{amsmath,amssymb,amsfonts}
\usepackage{algorithmic}
\usepackage{graphicx}
\usepackage{textcomp}

\newtheorem{remark}{Remark}
\newtheorem{definition}{Definition}

\begin{document}
\title{Sparse Structure Design for Stochastic Linear Systems via a Linear Matrix Inequality Approach}

\author{Yi Guo, 
Ognjen Stanojev, 
Gabriela Hug 
and Tyler Holt Summers
\thanks{This work is partially supported by an ETH Z\"{u}rich Postdoctoral Fellowship. This material is also based on work supported by the United States Air Force Office of Scientific Research under award number FA2386-19-1-4073 and by the National Science Foundation under award number ECCS-2047040 and CMMI-1728605. \textit{(Corresponding author: Yi Guo)}}
\thanks{Y. Guo, O. Stanojev and G. Hug are with Power Systems Laboratory at ETH Z\"{u}rich, Z\"{u}rich, 8092, Switzerland, email:\{guo, stanojev, hug\}@eeh.ee.ethz.ch.}
\thanks{T.H. Summers is with the Department
of Mechanical Engineering, The University of Texas at Dallas, Richardson,
TX, USA, email: tyler.summers@utdallas.edu. }
}

\maketitle

\begin{abstract}
In this paper, we propose a sparsity-promoting feedback control design for stochastic linear systems with multiplicative noise. The objective is to identify a sparse control architecture that optimizes the closed-loop performance while stabilizing the system in the mean-square sense. The proposed approach approximates the nonconvex combinatorial optimization problem by minimizing various matrix norms subject to the Linear Matrix Inequality (LMI) stability condition. We present two design problems to reduce the number of actuators via the static state-feedback and a low-dimensional output. A regularized linear quadratic regulator with multiplicative noise (LQRm) optimal control problem and its convex relaxation are presented to demonstrate the tradeoff between the suboptimal closed-loop performance and the sparsity degree of control structure. Case studies on power grids for wide-area frequency control show that the proposed sparsity-promoting control can considerably reduce the number of actuators without significant loss in system performance. The sparse control architecture is robust to substantial system-level disturbances while achieving mean-square stability.
\end{abstract}


\section{Introduction}\label{sec:intro}
\IEEEPARstart{D}{ynamical} systems with multiplicative noise provide rich models for many practical applications, including frequency control of power grids, deployment of robot agent teams, control of segmented mirrors in extremely large telescopes, and other applications in biological movement systems and aerospace engineering systems \cite{harris1998signal,todorov2005stochastic,todorov2003optimal,todorov2002optimal}. Especially in large-scale systems, substantial system-level disturbances and uncertainties may lead to oscillations and possibly instability. Therefore, one of the major challenges is to design efficient, high-performance, and robust control architectures that limit the number of actuators, sensors and actuator-sensor communication links to reduce the complexity and cost. 



Several recent efforts have demonstrated that sparse controller architectures can successfully balance closed-loop performance and controller complexity \cite{norris1989selection,fardad2011sparsity,schuler2011design,matni2017communication,fardad2014design,matni2016regularization,taha2018time,summers2016actuator,jovanovic2016controller,zare2018proximal,gravell2019sparse,summers2015submodularity,chang2018co,kim1991measure} and are therefore crucial to obtain efficient controllers for emerging complex systems. However, the design of these control architectures requires solving mixed combinatorial optimization problems.
One line of research formulated convex structured optimal control problems for controller design, such as symmetric modifications \cite{fardad2011sparsity}, diagonal modifications for optimal sensor and actuator selection \cite{dhingra2014admm} and a linear matrix inequality (LMI) approach with $\ell_1$-optimization \cite{polyak2013lmi}. Another line of research employed an algorithmic approach to solve the convex problems, such as the alternating direction method of multipliers \cite{dhingra2016method}, the proximal gradient and Newton methods \cite{moghaddam2016customized}, and also the second order method of multipliers for efficiently identifying the controller structure and its structured feedback synthesis \cite{dhingra2016method}. However, none of these works consider the multiplicative noise, which normally capture the disturbances on system model and the inherent uncertainties within input-output communication channels. System-level disturbances inherently appear on the system parameters and have fundamentally different effects on the state evolution than additive noise. In particular, a noise-ignorant classical optimal linear-quadratic controller may destabilize a stochastic system with multiplicative noise in the mean-square sense.
 
Linear systems with multiplicative noises are particular attractive as a stochastic modeling framework because they remain simple enough to admit closed-form expressions for stabilization \cite{boyd1994linear} and optimal control \cite{wonham1967optimal,kleinman1969optimal}. Many recent works proposed various approaches for controlling and filtering for systems with multiplicative noise, including an LMI approach \cite{gershon2001h}, the Riccati difference equation method \cite{rami2001solvability} and a game theory approach \cite{phillis1985controller}. However, all of this work uses fully populated control architectures, which become impractical and expensive as scale and complexity increase. These limitations of fully populated architectures motivate us to derive an approach that provide a sparse control architecture design for stochastic linear systems with multiplicative noise. 

Instead of performing a computationally-expensive combinatorial search, our approach leverages the convexity of various sparsity-promoting matrix norms to encourage of the feedback control matrix, while stabilizing the systems via LMI constraints. We first present two sparsity exploration problems to reduce the number of actuators by promoting static state-feedback and output-feedback with a low-dimensional output, with the goal to stabilize the stochastic system with multiplicative noise. We then formulate a regularized linear quadratic regulator with multiplicative noise (LQRm) optimal control problem and write it as a semi-definite programming (SDP) problem. This formulation tradeoffs the system performance and sparsity degree of the control architecture by different sparsity measures. Finally, we apply our approach to design a sparse wide-area frequency control structure in power grids. The numerical results on a four-bus system show that the control structure can be sparse at the expense of a slight loss in the performance. We also test the computational performance of our approach on an IEEE 39-bus network to visualize the tradeoffs under various levels of multiplicative noise. The noise-aware sparse structure requires more actuators to stabilize the system in the mean-square sense than the noise-ignorant design, which emphasizes the necessity of having an optimal structure approach for dynamical systems with significant system-level disturbances.

The rest of the paper is organized as follows: Section II presents our sparsity-promoting structure design approach for stochastic linear systems in the continuous-time domain. The objective is to stabilize the system with a reduced number of actuators via state-feedback control and output-feedback control with a low-dimensional output. A regularized LQRm problem and its relaxation are provided. Section III presents some case studies on power networks, while Section IV concludes the paper. Finally, the corresponding formulations for discrete-time stochastic linear systems are provided in the Appendix.

\textit{Notation}: We use $\mathbb{R}$, $\mathbb{R}_+$ and $\mathbb{R}_{++}$ to denote the sets of real numbers, non-negative real numbers and positive real numbers, respectively. Sets $\mathbb{S}^n$, $\mathbb{S}_+^n$ and $\mathbb{S}_{++}^n$ collect all $n$-dimensional symmetric matrices, semi-definite positive matrices and positive definite matrices, respectively. Given a matrix $M$, $M^\top$ denotes its transpose and $\textbf{Tr}(M)$ denotes its trace. We write $M \succeq 0$ ($M \succ 0$) to denote that $M$ is semi-positive definite (positive definite). For a given column vector $x\in\mathbb{R}^n$, we definite $\|x\|_1 := \sum_{i} |x_i|$, $\|x\|_2 := \sqrt{x^\top x}$ and $\|x\|_\infty := \max_i x_i$. Further, $|\cdot|$ denotes the absolute value of a number or the cardinality of a set, $\textbf{diag}(\cdot)$ constructs a diagonal matrix from a vector and
$\textbf{blkdiag}(\cdot)$ returns a block diagonal matrix. Lastly, $I$ denotes the identity matrix of appropriate dimension.

\section{Problem Formulation}\label{sec:formulations}

\subsection{Stochastic Linear Systems with Multiplicative Noise}
Consider a stochastic linear system with state- and input-dependent multiplicative noises \cite{arnold1974stochastic}:
\begin{equation}\label{eq:ito_equation}  
    dx_t = (A_0x_t + B_0u_t) dt + \sum_{i=1}^k\sigma_iA_i x_t d\beta_{it} + \sum_{j = 1}^\ell\rho_j
    B_j u_t d\delta_{jt},
\end{equation}
where $x_t \in \mathbb{R}^n$ denotes the state vector, $u_t\in\mathbb{R}^m$ denotes the control input vector and $\beta_{it}(i = 1,\ldots,k)$ and $\delta_{jt}(j = 1,\ldots,\ell)$ denote the disturbances. We assume these disturbances to be zero mean uncorrelated stationary normalized Wiener processes. The following properties hold:
\begin{equation*}
\begin{aligned}
    \mathbf{E}[d\beta_{it}] = 0,~~~ \mathbf{E}[d\delta_{jt}] = 0,~~~
    \mathbf{E}[d\beta_{it}^2] = dt,~~~
    \mathbf{E}[d\delta_{jt}^2] = dt,\\
    \mathbf{E}[d\beta_{i_1t}d\beta_{i_2t}] = 0~ (i_1 \neq i_2),~  \mathbf{E}[d\delta_{j_1t}d\delta_{j_2t}] = 0~ (j_1 \neq j_2),\\
    \textrm{and} ~~~~~~~~~~~~~~~~~~~~~~\mathbf{E}[d\beta_{it}d\delta_{jt}] = 0,~~~~~~~~~~~~~~~~~~~~~~~~~\\
    \forall i,i_1,i_2 = 1,\ldots,k, ~\textrm{and} ~~\forall j,j_1,j_2 = 1,\ldots,\ell.~~~~~~~~
\end{aligned}
\end{equation*}
\noindent The scale factors $\sigma_i$ and $\rho_j$ indicate the intensities of the disturbances, which scale the unit variance of $d\beta_{it}$ and $d\delta_{jt}$. The initial condition of the system \eqref{eq:ito_equation} is $\mathbf{E}[x_0] =  0$ and $\mathbf{E}[x_0 x_0^\top] =  \Sigma_0$, where $\Sigma_0\in\mathbb{R}^{n\times n}$ is the covariance matrix of the initial state vector. The matrices $A_0 \in \mathbb{R}^{n\times n}$ and $B_0\in\mathbb{R}^{n \times m}$ correspond to the constant system matrices. The state diffusion term projects state-dependent noise $d\beta_{it}$ by matrix $A_i\in \mathbb{R}^{n \times n}$, and the input diffusion term projects input-dependent noise $d\delta_{jt}$ by matrix $B_j \in \mathbb{R}^{n \times m}$. Assuming the dynamic system \eqref{eq:ito_equation} is open-loop mean-square unstable, we apply the sparse ideology in the design of optimal linear feedback control while stabilizing the system in the mean-square sense. We first present a LMI condition for the mean-square stability of \eqref{eq:ito_equation} and then discuss how to find out a subset of actuators to trade off the system closed-loop performance with various degrees of sparsity.


\begin{definition}
The system \eqref{eq:ito_equation} is \emph{mean-square stable} if for every initial condition $\mathbf{E}[x_0 x_0^\top] = \Sigma_0$ the solution of \eqref{eq:ito_equation} satisfies
\begin{equation*}
    \lim_{t \to +\infty}\mathbf{E}[x_t^\top x_t] = 0.
\end{equation*}
\end{definition}
Note that the above mean-square stability condition will converge to a constant if the system \eqref{eq:ito_equation} also have non-zero mean additive noise.

\begin{definition}
The system \eqref{eq:ito_equation} with an initial condition $\mathbf{E}[x_0x_0^\top] = \Sigma_0$ is called (mean-square) stabilizable if there exists a mean-square stabilizing state-feedback control in the form $u_t = Kx_t$,  where $K$ is a constant matrix.
\end{definition}

\begin{definition}
Given the output of the system \eqref{eq:ito_equation},  such that $y_t = Cx_t$, the system \eqref{eq:ito_equation} with an initial condition $\mathbf{E}[x_0x_0^\top] = \Sigma_0$ is called (mean-square) stabilizable if there exists a mean-square stabilizing output-feedback control in the form $u_t = Ky_t$, where $K$ is a constant matrix.
\end{definition}

\subsection{Stabilization with a Reduced Number of State-Feedback Controllers}
Assume the system \eqref{eq:ito_equation} is open-loop mean-square unstable and stabilizable via the state-feedback control, the goal of this subsection is to identify potential row-sparsity patterns of the closed-loop state-feedback control law in the form of $u_t = Kx_t$, such that the closed-loop system described by
\begin{equation}\label{eq:ito_equation_closed_loop}  
    dx_t = (A_0 + B_0K)x_t dt + \sum_{i=1}^k\sigma_iA_i x_t d\beta_{it} + \sum_{j = 1}^\ell\rho_{i}
    B_j Kx_t d\delta_{jt},
\end{equation}
is mean-square stable. To guarantee the mean-square stability of the closed-loop system \eqref{eq:ito_equation_closed_loop}, the static state-feedback control gain matrix $K$ exists if and only if there exists a matrix $X\in\mathbb{S}_{++}^n$ such that the following condition holds \cite{el1995state}:
\begin{equation}\label{eq:LMI_lyapunov_orginal}\nonumber
\begin{aligned}
    (A_0 + B_0K)^\top X + X(A_0 + B_0K) + \sum_{i = 1}^k \sigma_i^2A_i^\top X A_i \\
    + \sum_{j = 1}^\ell \rho_j^2 K^\top B_j^\top X B_jK \prec 0.
\end{aligned}
\end{equation}
Pre- and post-multiplying the above inequality by $P = X^{-1}$ and introducing a new variable $Y = KP$, we arrive at the following condition:
\begin{equation}\label{eq:continuous_stable}
\begin{aligned}
    A_0P + PA_0^\top + B_0Y + Y^\top B_0^\top + \sum_{i = 1}^k \sigma_i^2 (A_iP)^\top P^{-1} A_iP \\
    + \sum_{j = 1}^\ell \rho_j^2  (B_jY)^\top P^{-1} B_j Y \prec 0,
\end{aligned}
\end{equation}
where $P = P^\top \in\mathbb{S}^n_{++}$ and $Y\in \mathbb{R}^{m\times n}$ are the matrix variables. A stabilizing state-feedback controller can be reconstructed by $K = YP^{-1}$. Leveraging the Schur's Lemma, we transform the condition \eqref{eq:continuous_stable} together with $P \succ 0$ into a LMI:
\begin{equation} \label{eq:LMOstabilityContinuous}
    \begin{bmatrix}
    A_0P + PA_0^\top + B_0Y + Y^\top B_0^\top & Z\\
    Z^\top & Z_P
    \end{bmatrix} \prec 0,
\end{equation}
where $Z = \begin{bmatrix} \sigma_1 PA_1^\top, \ldots, \sigma_k PA_k^\top, \rho_1 Y^\top B_1^\top,\ldots,\rho_\ell Y^\top B_\ell^\top \end{bmatrix}$ and 
\begin{equation*}
Z_P = \textbf{blkdiag}\underbrace{\left(-P, \ldots, -P\right)}_{k+\ell}.
\end{equation*}
If $Y$ is row sparse, then the state-feedback law $K$ is row sparse as well since post-multiplication preserves the zero-row structure. Hence, we promote row sparsity of $Y$ through the following SDP:
\begin{equation}\label{eq:rsp_problem}
    \min_{Y,P} ~\|Y\|_r, \quad \textrm{subject to:}\quad \eqref{eq:LMOstabilityContinuous} \quad \textrm{and} \quad P \succ 0,
\end{equation}
where $\|Y\|_r$ represents a generic row-sparsity induced function that can be chosen from row-norm \cite{polyak2013lmi}, group LASSO\cite{ahsen2017two} and sparse Group LASSO \cite{ahsen2017two}. Note that the stabilizing state-feedback control matrix $K_{\textrm{rsp}}$ with an identified row sparse pattern can be obtained from the solution $P^*_\textrm{rsp}$, $Y^*_\textrm{rsp}$ of the above SDP problem, with the linear feedback control matrix calculated as $K_{\textrm{rsp}} = Y^*_\textrm{rsp}{P^*_\textrm{rsp}}^{-1}$.

\begin{remark}[Sparsity-promoting Norms]We now suggest different sparsity-promoting norms that can be used in \eqref{eq:rsp_problem} and the problems will be presented in the following section. 
\begin{itemize}
\item \textit{\underline{Row Norm and Column Norm} \cite{polyak2013lmi}:} Given a matrix $Y\in\mathbb{R}^{m\times n}$, the row and column sparsity can be induced by various sparsity-promoting matrix norms respectively defined as:
\begin{equation*}
    \|Y\|_{\textrm{row}} = \sum_{i=1}^m \|Y_{r,i}\|_{\infty}, \quad \|Y\|_{\textrm{col}} = \sum_{j=1}^n \|Y_{c,j}\|_{\infty},
\end{equation*}
where $\|Y_{r,i}\|_{\infty}$ and $\|Y_{c,j}\|_{\infty}$ are the maximum absolute values of the $i$-th row and $j$-th column of matrix $Y$, respectively. 

\item \textit{\underline{Group LASSO}\cite{ahsen2017two}:} Row and column sparsity can also be induced by the row and column group LASSO:
\begin{equation}\nonumber
        \|Y\|_{\textrm{rGL}} = \sum_{i=1}^m \|Y_{r,i}\|_2,  \quad \|Y\|_{\textrm{cGL}} = \sum_{j=1}^n \|Y_{c,j}\|_2,
\end{equation}
where $\|Y_{r,i}\|_2$  and $\|Y_{c,j}\|_2$ are the vector $\ell_2$-norms of the $i$-th row and $j$-th column of matrix $Y$, respectively. 

\item\textit{\underline{Sparse Group LASSO}\cite{ahsen2017two}:} The row and column sparse group LASSO  can also promote the sparsity pattern:    
\begin{equation}\nonumber
\begin{aligned}
        &\|Y\|_{\textrm{rSGL},\mu} = \sum_{i=1}^m (1-\mu)\|Y_{r,i}\|_1 + \mu \|Y_{r,i}\|_2,\\
        &\|Y_c\|_{\textrm{cSGL},\mu} = \sum_{j=1}^n (1-\mu)\|Y_{c,j}\|_1 + \mu \|Y_{c,j}\|_2,
\end{aligned}
\end{equation}
where $\|Y_{r,i}\|_1$ and $\|Y_{c,j}\|_1$ are the vector $\ell_1$-norms of the $i$-th row and $j$-th column of matrix $Y$, respectively. The constant $\mu\in[0,1]$ quantifies the weight on the two combined norms. In the rest of this paper, we refer to $\|Y\|_{\textrm{reg}}$ as a generic sparsity-promoting regularizer in the following optimal design formulation.
\end{itemize}
\end{remark}

\subsection{Stabilization with a Reduced Number of Output-Feedback Controllers via a Low-Dimensional Output}
In this subsection, we present a stabilization solution to reduce the number of output-feedback controllers $u_t = Ky_t$ via a low-dimensional output $y_t = Cx_t$, where $y_t \in\mathbb{R}^{n_y}$ is the output vector and $C\in\mathbb{R}^{n_y\times n}$ is the output matrix. Assume that the system \eqref{eq:ito_equation} has exact measurements of the full states and is stabilizable via an output-feedback control law $K$. Note that the potential stabilization solutions via output-feedback controllers highly depend on the $C$ matrix. Here we assume that there is at least one low-dimensional output that can enable the sparse control structure.
The goal of this subsection is therefore to obtain a low-dimensional system output matrix $C$ and a column-row sparse output-feedback matrix $K$ to stabilize \eqref{eq:ito_equation}. To achieve this, we first change the row sparsity promoting function $\|Y\|_r$ in \eqref{eq:rsp_problem} to a generic column sparsity induced norm $\|Y\|_c$ resulting in
\begin{equation}\label{eq:csp_problem}
    \min_{Y,P} ~\|Y\|_c, \quad \textrm{subject to:}\quad \eqref{eq:LMOstabilityContinuous} \quad \textrm{and} \quad P \succ 0.
\end{equation}
Similarly to the row sparse state-feedback law from \eqref{eq:rsp_problem}, the solution of \eqref{eq:csp_problem} promotes the column sparsity on $Y_c^*$. Similar to the previous section, the feedback law is in the form:
\begin{equation}\label{eq:csp_state_feedback}
u_t = Y^*_\textrm{csp}{P^*_\textrm{csp}}^{-1}x_t, 
\end{equation}
where $Y^*_\textrm{csp}$ and $P^*_\textrm{csp}$ are now the solution of \eqref{eq:csp_problem}. Interestingly, the sparsity pattern of the output-feedback $u_t = K_{\textrm{csp}}y_t$ can be attained by mapping the matrix multiplication of the term $Y^*_\textrm{csp}{P^*_\textrm{csp}}^{-1}$ in \eqref{eq:csp_state_feedback} and the term ${K}_{\textrm{csp}}{C}_{\textrm{rsp}}$ in 
\begin{equation}\nonumber
    u_t = {K}_{\textrm{csp}}y_t = {K}_{\textrm{csp}}{C}_{\textrm{rsp}}x_t. 
\end{equation}
The output-feedback law $K_{\textrm{csp}}$ consists of the nonzero-columns of $Y_\textrm{csp}^*$ and the output matrix ${C}_{\textrm{rsp}}$ is composed of the rows of ${P^*_\textrm{csp}}^{-1}$ with the same indices. In this way, we reduce the number of necessary outputs in $y_t = {C}_{\textrm{rsp}}x_t$ while stabilizing the system by a column sparse output-feedback $u_t = {K}_{\textrm{csp}}y_t$. 

After the identification of the column sparsity pattern of the output-feedback law ${K}_{\textrm{csp}}$, we now explore the potential row sparsity to reduce the number of the output-feedback controllers. Having the knowledge of the column sparsity pattern of variable $Y$ by solving \eqref{eq:csp_problem}, we integrate the zero-column pattern of the solution $Y_\textrm{csp}^*$ as additional constraints into \eqref{eq:rsp_problem} resulting in:
\begin{subequations}\label{eq:rsp_problem_constraints}
\begin{align}
    & \min_{Y,P} && \|Y\|_r, \\
    & \textrm{subject to:} &&  Y_{c,i} = 0, ~\forall i \in \mathcal{C}, \label{eq:rsp_zero_column_constraints} \\
    & &&  P \succ 0 \quad \textrm{and} \quad \eqref{eq:LMOstabilityContinuous},
    \end{align}
\end{subequations}
where $Y_{c,i}$ indicates the $i$-th column of the variable $Y$ and the set $\mathcal{C}$ collects the indices of all zero columns of $Y_\textrm{csp}^*$. The solution of \eqref{eq:rsp_problem_constraints} are $Y_{\textrm{sp}}^*$ and ${P^*_\textrm{sp}}$. We adopt the output-feedback law ${K}_{\textrm{sp}}$ as the same row-column sparsity of $Y_{\textrm{sp}}^*$ and the output matrix ${C}_{\textrm{sp}}$ consists of the non-zero rows of ${P^*_\textrm{sp}}^{-1}$. The feedback control structure $u_t = K_{\textrm{sp}}y_t$ has the row-sparsity such that the corresponding controllers can be removed with a low-dimensional output $y_t= {C}_{\textrm{sp}}x_t$. Hence, by solving these two sparsity promoting problems \eqref{eq:csp_problem} and \eqref{eq:rsp_problem_constraints} in sequence, we can identify a sparsity structure to design low-dimensional outputs and remove output-feedback controllers that are not necessary for stabilizing the system. Clearly, the sparse control structure with only a subset of outputs and controllers will reduce system performance. We will discuss the performance degradation and sparsity tradeoffs later in the paper. It is worth emphasizing that the above sparse structure design problems are also applicable for open-loop mean-square stable systems with the goal to allow for performance tradeoffs, which will be discussed later in this paper.

\subsection{Tradeoffs between Optimal Control and Sparsity}
We now consider an application of our approach to the linear quadratic regulator with multiplicative noise (LQRm) for the system \eqref{eq:ito_equation} given an initial condition $\mathbf{E}[x_0x_0^\top] = \Sigma_0$:
\begin{subequations}\label{eq:LQRm}
\begin{align}
  &\min_{u_t(\cdot)}\quad  J(\Sigma_0, u(\cdot))  = \mathbf{E}\int_0^\infty \left(x^\top_t Q x_t + u^\top_t R u_t\right)dt,\\
  & \textrm{subject to}:\\
  & dx_t = (A_0x_t + B_0u_t) dt + \sum_{i=1}^k\sigma_iA_i x_t d\beta_{it} + \sum_j^\ell\rho_j
    B_j u_t d\gamma_{jt}, \label{eq:LQRm_dynamics}
\end{align}
\end{subequations}
where $Q$ and $R$ are positive definite. The objective is to determine an optimal linear state-feedback (output-feedback) control law that trades off the LQRm closed-loop performance and the sparsity of the linear feedback control gain $K$\footnote{We initially describe the regularized LQRm problem with static state-feedback control. The corresponding reformulation with static output-feedback control follows closely the line of the steps as state-feedback control but having the output-feedback control law $u_t = Ky_t = KCx_t$.}. For state-feedback control $u_t = Kx_t$, we are ultimately interested in a regularized LQRm problem with an alternative objective:
\begin{equation}\nonumber
    \min_{K} \quad  J(\Sigma_0, u_t(K)) + \gamma J_{\textrm{reg}}(K),
\end{equation}
where $J_{\textrm{reg}}(K)$ is a sparsity-promoting function of $K$ and $\gamma$ specifies the importance of its sparsity. Together with the stability constraint \eqref{eq:LMOstabilityContinuous}, the optimal control problem of determining the stabilizing closed-loop state-feedback $K$ that minimizes the LQRm cost and determines the potential sparsity structure of $K$ can be reformulated as the following optimization problem \cite{boyd1994linear}:
\begin{subequations}\label{eq:LQR_nontractable}
    \begin{align}
    & \min_{Y,P} \quad \quad \mathbf{Tr}(P^{-1}\Sigma_0) + \gamma \|Y\|_\textrm{reg}, \label{eq:LQR_nontractable_obj}\\
    & \textrm{subject to:} \quad P\succ 0,\\
    & A_0P + PA_0^\top + B_0Y + YB_0^\top + \sum_{i=1}^k \sigma_i PA_i^\top P^{-1} A_i P \nonumber \\
    & + \sum_{j=1}^l \beta_j Y^\top B_j^\top P^{-1} B_j Y + Y^\top R Y + PQP \prec 0, \label{eq:LQR_nontractable_constraints}
    \end{align}
\end{subequations}
Note that \eqref{eq:LQR_nontractable} is intractable due to the matrix inverse $P^{-1}$ in the objective and the nonlinear multiplication in \eqref{eq:LQR_nontractable_constraints}. We now introduce new (slack) variables $\Pi \in \mathbb{R}^{n\times n}$ and $\kappa \in \mathbb{R}$ and introduce a constraint that provides an upper bound of the LQRm cost, i.e., $J = \mathbf{Tr}(\Sigma_0P^{-1})\leq \kappa$. By leveraging the Schur's Lemma \cite{boyd1994linear,rami2000linear}, we obtain to the following SDP problem:
\begin{subequations}\label{eq:LQR_SDP}
    \begin{align}
    & \min_{Y,P,\Pi,\kappa} \quad \quad \kappa + \gamma \|Y\|_\textrm{reg},\\
    & \textrm{subject to:} \quad\quad \mathbf{Tr}(\Pi) \leq \kappa, \begin{bmatrix}
    \Pi & \Sigma^{\frac{1}{2}}_0\\
    \Sigma^{\frac{1}{2}}_0 & P
    \end{bmatrix} \succeq 0, P \succ 0,\\
    & \scalebox{0.94}[1]{$\begin{bmatrix} A_0P + PA_0^\top + B_0Y + Y^\top B_0^\top & Z & Y^\top & P \\
    Z^\top & Z_P & 0 & 0\\
    Y & 0 & - R^{-1} & 0\\
    P & 0 & 0 & - Q^{-1}
    \end{bmatrix} \prec 0.$}
    \end{align}
\end{subequations}
The solution ${Y}^*$, ${P}^*$, ${\Pi}^*$ and ${\kappa}^*$ defines a sparse stabilizing law ${K}_{\textrm{sp}} = {Y}^*{P^*}^{-1}$ and the upper bound of the LQRm cost $J^*(K) = \mathbf{Tr}(\Sigma_0{P^*}^{-1})\leq {\kappa}^*$. Note that $\eqref{eq:LQR_SDP}$ is convex and can be solved by several academic and commercial SDP solvers. The proposed method identifies a sparse control structure with only a small performance loss, as demonstrated by the numerical studies in the next section.

\section{An Application to Power Grids}\label{sec:case_studies}
We apply the proposed methodology to devise an optimal wide-area frequency control scheme for a power transmission system. The objective is to design a sparse linear feedback frequency control architecture, which stabilizes the frequency dynamics with modeling errors as multiplicative noise. 
\subsection{Network Modeling}
Consider a lossless transmission system modeled as a graph $G = (\mathcal{N},\mathcal{E})$ with nodes (or buses) $\mathcal{N} = \{1,\cdots, N\}$ and edges (or lines) $\mathcal{E} \subseteq \mathcal{N}\times \mathcal{N}$. The topology of the grid is represented by the Laplacian matrix $L\in\mathbb{R}^{N\times N}$ induced by the line susceptances $b_{ij}$ for all $(i,j)\in\mathcal{E}$ (see \cite{guo2019performance}). We partition all buses into a set of buses with generators $\mathcal{G}$ (i.e., synchronous machines and inverter-based generators) and a set of buses with frequency-sensitive loads $\mathcal{L}$, where $\mathcal{N} = \mathcal{G}\cup \mathcal{L}$. In this paper, we consider a reduced-order network model such that we merge the buses only having frequency-insensitive loads or having no load with the frequency-sensitive buses. The variables describing the state of the network include angles $\theta_i$ and frequency $\omega_i$ for $i\in\mathcal{N}$. The associated system dynamics \cite{hill1990stability,machowski2020power} derived from the linearized swing equations are given by 
\begin{subequations} \label{eq:DAE_power_grids}
\begin{align}
\dot{\theta}_i & = \omega_i, &&\forall i \in \mathcal{N},\\
M_i \ddot{\theta}_i + D_{g,i}\dot{\theta}_i & = - \sum_{(i,j)\in\mathcal{E}} b_{ij}(\theta_i - \theta_j) + u_i, \hspace{-3mm}&&\forall i\in\mathcal{G},  \label{eq:DAE_power_grids_G}\\
D_{l,i} \dot{\theta}_i & =  - \sum_{(i,j)\in\mathcal{E}} b_{ij}(\theta_i - \theta_j) + u_i, \hspace{-3mm}&&\forall i \in \mathcal{L}, 
\end{align}
\end{subequations}
where $u_i\in\mathbb{R}$ is a controllable generation or load for all $i\in\mathcal{N}$. A generator bus $i\in\mathcal{G}$ is characterized by its inertia $M_i \in\mathbb{R}_{++}$ (rotational or virtual inertia) and its droop coefficient $D_{g,i} \in\mathbb{R}_{++}$. A frequency-sensitive load bus $i\in\mathcal{L}$ is characterized by its sensitivity coefficient $D_{l,i} \in\mathbb{R}_{++}$. 

Here, we consider inertia variations caused by inverter-based generation or modeling errors, which are modeled by treating the inertia parameters $M_i$ as multiplicative noise rather than simply a constant. The inertia parameters for all $i\in\mathcal{G},$ $M_i$ are therefore modeled as random parameters with the mean value $\overline{M}_i$ and the variance $\overline{\sigma}_i^2$ \cite{guo2019performance}. In general, the modeling errors of the network topology $L$ and frequency-sensitive coefficients $D_{\ell,i}/D_{g,i}$ can also be treated as multiplicative noises in \eqref{eq:DAE_power_grids}.  However, for simplicity, only the inertia parameter randomness is considered. The above dynamics can be formulated in a generalized as a multi-input multi-output stochastic linear system with multiplicative noise:
\begin{equation}\label{LTVstateform_reduced}
\begin{aligned}
& dx_t = A_0x_tdt + \sum_{i = 1}^{|\mathcal{G}|} \sigma_i \left(A_i x_t +  B_i u_t\right)d\beta_{it},\\
& y_t = Cx_t, ~~~~~~~x_t = 
\begin{bmatrix}
\theta^\top_{g,t}, \omega^\top_{g,t}, \theta_{l,t}^\top
\end{bmatrix}^\top,
\end{aligned}
\end{equation}
where
\begin{equation*} 
\scalebox{0.95}[1]{$
\begin{aligned}
& A_0  = 
\begin{bmatrix}
0 & I & 0\\
-M^{-1}L_{gg} & -M^{-1}D_g & -M^{-1}L_{gl}\\
-D^{-1}_lL_{l g} & 0 & -D^{-1}_lL_{ll}
\end{bmatrix},\\ 
& B_0 = \begin{bmatrix}
  0 \\
  M^{-1}\\
  0
\end{bmatrix}, A_i =
\begin{bmatrix}
0 & 0 & 0\\
R_iL_{gg} & R_iD_g & R_iL_{g\ell}\\
0 & 0 & 0
\end{bmatrix}, B_i = \begin{bmatrix}
  0 \\
  R_i\\
  0
\end{bmatrix},
\end{aligned}
$}
\end{equation*}
and the vectors $\theta_{g,t}\in\mathbb{R}^{|\mathcal{G}|}$/$\theta_{l,t}\in\mathbb{R}^{|\mathcal{L}|}$ and $\omega_{g,t}\in\mathbb{R}^{|\mathcal{G}|}$ collect the angle states of generation/load buses and the frequency states of generation buses, respectively. Note that the inertia parameter for every generator $M_i$, $i\in\mathcal{G}$ appears in the state compact form \eqref{LTVstateform_reduced} as the inverse distribution of $M_i$ with the mean value $\overline{M}^{-1}_i$ and the variance $\sigma_i^2$. The matrix $M^{-1} := \textbf{diag}(\overline{M}_i^{-1})\in\mathbb{R}^{|\mathcal{G}|\times|\mathcal{G}|}$ collects the mean values of the inverse distribution of $M_i$. The modeling error of the inverse of inertia parameter $d\beta_t$ for all $i\in\mathcal{G}$ is considered as an independent Wiener process normalized by $\sigma_i$. The diagonal matrix $D_g \in \mathbb{R}^{|\mathcal{G}|\times|\mathcal{G}|}$ collects the droop coefficients $D_{g,i}$ at all buses $i\in\mathcal{G}$ and the diagonal matrix $D_l \in \mathbb{R}^{|\mathcal{L}|\times|\mathcal{L}|}$ collects the sensitivity coefficients $D_{l,i}$ at all buses $i\in\mathcal{L}$. The matrices $L_{gg}\in\mathbb{R}^{|\mathcal{G}|\times|\mathcal{G}|}$, $L_{lg}\in\mathbb{R}^{|\mathcal{L}|\times|\mathcal{G}|}$ and $L_{gl}\in\mathbb{R}^{|\mathcal{G}|\times|\mathcal{L}|}$, $L_{ll} \in\mathbb{R}^{|\mathcal{L}|\times|\mathcal{L}|}$ are derived from the original Laplacian matrix $L$, which represent the weighted connections between different types of buses (i.e., load and generation). We define an \emph{inertia disturbance allocation} matrix $R_i \in \mathbb{R}^{N\times N}$ in $A_i$ and $B_i$ associated with each bus $i\in\mathcal{G}$. The elements in $R_i$ are all zeros except for one diagonal element $r_{ii} = 1$, which maps the corresponding inertia disturbance $d\beta_i$ onto bus $i$. If the inertia variation at bus $i$ is insignificant, we set $R_i = 0$ to remove the inertia disturbance on the $i$-th bus. 

We conduct the numerical experiments on two power networks. A small, four-bus power system is used to numerically show the design results and a larger scale power system is utilized to demonstrate the computational efficiency of the proposed LMI approach. We start with the four-bus system for which the grid topology and line parameters can be found in \cite{guo2019performance}. The damping/frequency sensitivity coefficients are set to $D_i = 10$ for all buses $i\in\mathcal{N}$. The mean value of inertia is $\overline{M}_i = 10$ for all $i\in\mathcal{G}$ and the standard deviation of the distribution of $M_i^{-1}$ is 10\% of the mean value, i.e., $\sigma_i = 10\% \overline{M}_i^{-1}$ for all $i\in\mathcal{G}$. Note that the open-loop dynamics \eqref{eq:DAE_power_grids} are not mean-square stable if the multiplicative noise variance is significant. More details with respect to the modeling and stability analysis related to this example can be found in \cite{guo2019performance}. In this paper, we mainly focus on finding the sparse structure of the closed-loop feedback control for generators that stabilize the frequency dynamics with inertia disturbances. The initial condition of the states is $\Sigma_0 = 0.1I$ and the coefficients of the LQRm cost are $Q = I$ and $R = I$.
\subsection{Reducing the Number of Actuators via Static State-Feedback}
We first check if the stochastic linear system \eqref{LTVstateform_reduced} is mean-square stabilizable (no sparsity induced) via a closed-loop state-feedback controller $u_t = Kx_t$ by solving \eqref{eq:LQR_SDP} with $\gamma = 0$. The obtained solution is
\begin{equation*}
\scalebox{0.89}[1]{$
\begin{aligned}
&K_0 =\\
&\begin{bmatrix}
   -0.2329 & -0.0939 & -0.0236 & -0.2740 & -0.0879 & -0.0216\\
   -0.0941 & -0.3008 & -0.0823 & -0.0878 & -0.3358 & -0.0770\\
   -0.0236 & -0.0822 & -0.2075 & -0.0216 & -0.0771 & -0.2503
    \end{bmatrix}.
\end{aligned}
$}
\end{equation*}
With bus 4 grounded as an infinite bus and employing model reduction \cite{guo2019performance}, the system has six states, three inputs and three independent multiplicative noises and all three generators participate in stabilizing the grid. The state-feedback law $K$ is fully populated and the LQRm cost is $\kappa_0^* = 1.7915$. To obtain a row sparse solution $K$ and therefore a reduced number of state-feedback controllers, we again solve the regularized LQRm \eqref{eq:LQR_SDP} using the row-norm $\|Y\|_{\textrm{row}}$ with $\gamma = 4$, which results in the following sparse structure:
\begin{equation*}
\scalebox{0.89}[1]{$
\begin{aligned}
&{K}_{\textrm{rsp},4} =\\
&\begin{bmatrix}
    -0.0114 & -0.0153 & -0.0103 & -0.0110 & -0.0147 & -0.0099\\    
    -0.0718 & -0.0960 & -0.0645 & -0.0692 & -0.0925 & -0.0622\\
    0 & 0 & 0 & 0& 0 & 0
    \end{bmatrix}.
\end{aligned}
$}
\end{equation*}
This leads to a slight increase in the LQRm cost, namely $\kappa_4^* = 1.8757$. The row sparse state-feedback control law $K_{\textrm{rsp},4}$ indicates that generator 3 is not necessary to stabilize the system but at the expense of 4.67\% decrease of the closed-loop performance.
 
\subsection{Reducing the Number of Actuators via a Low-Dimensional Output}
To obtain a sparse structure of the output-feedback control $u_t = Ky_t = KCx_t$, we first solve \eqref{eq:LQR_SDP} by using the column-norm $\|Y\|_{\textrm{col}}$ as the sparsity regularizer with $\gamma = 0.5$. The solution is given by
\begin{equation*}
Y^*_{\textrm{csp},0.5} = 
\begin{bmatrix}
-0.0177 & -0.0400 & ~0.0078 & 0 & 0 & 0\\
-0.0177 & -0.0439 & -0.0080 & 0 & 0 & 0\\
~0.0074 & -0.0387 & -0.0080 & 0 & 0 & 0
\end{bmatrix}.
\end{equation*}
We adopt the sparse output-feedback law by taking the non-zero columns of $Y^*_{\textrm{csp},0.5}$:
\begin{equation*}
K_{\text{csp}, 0.5} =
\begin{bmatrix}
-0.0177 & -0.0400 & ~0.0078\\
-0.0177 & -0.0439 & -0.0080 \\
~0.0074 & -0.0387 & -0.0080
\end{bmatrix}.
\end{equation*}
The associated three-dimensional output matrix consists of the first three rows of the solution  ${P_\textrm{csp,0.5}^*}^{-1}$, namely:
\begin{equation*}
\scalebox{0.90}[1]{$
C_{\textrm{rsp},0.5} =
\begin{bmatrix}
    2.9985 & 0.9354 & 0.2549 & 2.3444 & 0.9222 & 0.2352\\
    0.9354 & 3.6956 & 0.8172 & 0.9155 & 3.0168 & 0.8046\\
    0.2549 & 0.8172 & 2.7508 & 0.2349 & 0.8127 & 2.0984
    \end{bmatrix}.
$}
\end{equation*}

We now reduce the number of output-feedback controllers by exploring the row sparsity of $K_{\textrm{csp}}$. Forcing the last three columns of $Y$ to be equal to zero as additional constraints, we solve \eqref{eq:LQR_SDP} again using the row-norm regularizer and $\gamma = 2.6$. The solution $Y^*_{\textrm{sp},2.6}$ has the following row-column sparse pattern:
\begin{equation*}
Y^*_{\textrm{sp},2.6} = 
\begin{bmatrix}
0 & 0 & 0 & 0 & 0 & 0\\
-0.0130 & -0.0130 & -0.0130 & 0 & 0 & 0\\
0 & 0 & 0 & 0 & 0 & 0
\end{bmatrix}.
\end{equation*}
Hence, we adopt
\begin{equation}\nonumber
    K_{\textrm{sp},2.6} = 
    \begin{bmatrix}
    0 & 0 & 0\\
    -0.0130 & -0.0130 & -0.0130\\
    0 & 0 & 0
    \end{bmatrix},
\end{equation}
as a row-column sparse output-feedback law. The output matrix $C_{\textrm{sp},2.6}$ is composed of the first three rows of ${P_{\textrm{sp},2.6}^*}^{-1}$, namely:
\begin{equation*}
\scalebox{0.90}[1]{$
C_{\textrm{sp},2.6} =
\begin{bmatrix}
    3.0963 & 1.0656 & 0.3010 & 2.4546 & 1.0592 & 0.2839\\
    1.0656 & 3.8739 & 0.9181 & 1.0653 & 3.2085 & 0.9221\\
    0.3010 & 0.9181 & 2.8063 & 0.2843 & 0.9176 & 2.1612
    \end{bmatrix}.
$}
\end{equation*}
At the end, the designed output-feedback controller uses one generator with a three-dimensional output feedback to stabilize the system at the expense of 5.5\% LQRm cost increase compared to $\kappa_0^* = 1.7915$.

\subsection{Tradeoff System Performance and Degree of Sparsity }
To discuss the computational cost of our approach on a large-scale system, we use the IEEE 39-bus New England transmission system to demonstrate the tradeoffs between the sparsity degree of the structure and the LQRm cost under various multiplicative noise settings. This system consists of 39 buses and 10 generators, where generator 10 is an equivalent aggregated model \cite{zimmerman2010matpower}. The inertia mean value is $M_i = 10$ for all generators $i\in\mathcal{G}$ and the damping/frequency sensitivity coefficient is set to $D_i = 10$ for all buses $i\in\mathcal{N}$. We vary the sparsity importance $\gamma$ to tradeoff the sparsity degree of the state-feedback control law $K$ and the LQRm performance. The LQRm cost coefficients are $Q = I$ and $R = I$. The initial state condition is $\Sigma_0 = 0.1I$. We used the MOSEK SDP solver \cite{aps2019mosek} via the MATLAB interface CVX \cite{grant2014cvx} on a laptop with 16 GB memory and 2.3 GHz Intel Core i7-10510U CPU. It took 121.7 seconds to solve \eqref{eq:LQR_SDP} with 10 inputs, 47 states and 10 independent multiplicative noises. 

\begin{figure}
\centering
\includegraphics[width=3.5in]{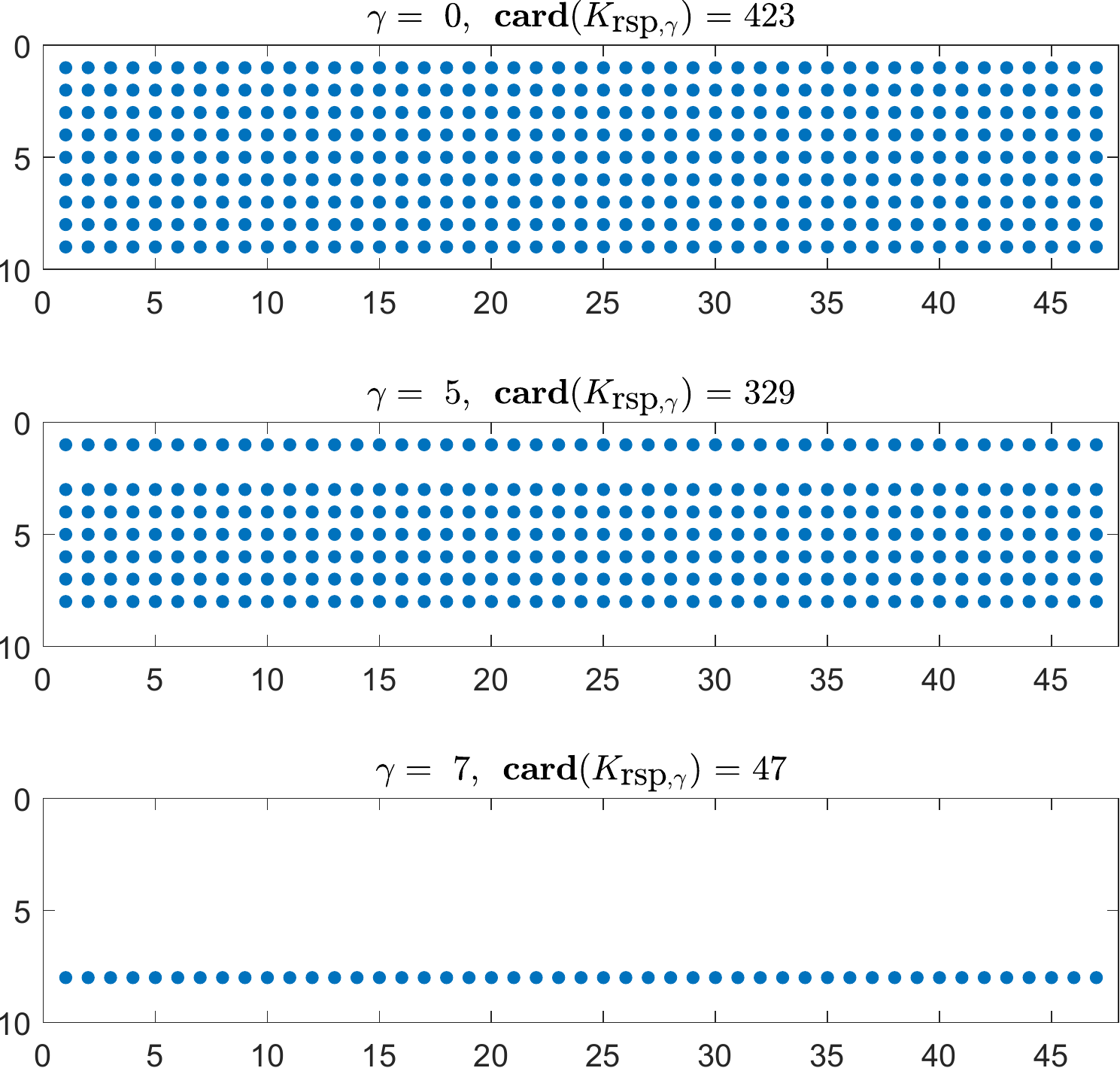}
\caption{Sparsity patterns of the closed-loop state-feedback $K$ resulting from the row-norm regularizer under various $\gamma$. The white elements represent (nearly) zero entries}
\label{fig:different_sparsity_K}
\end{figure}

\begin{figure}[htbp!]
\centering
\includegraphics[width=3.5in]{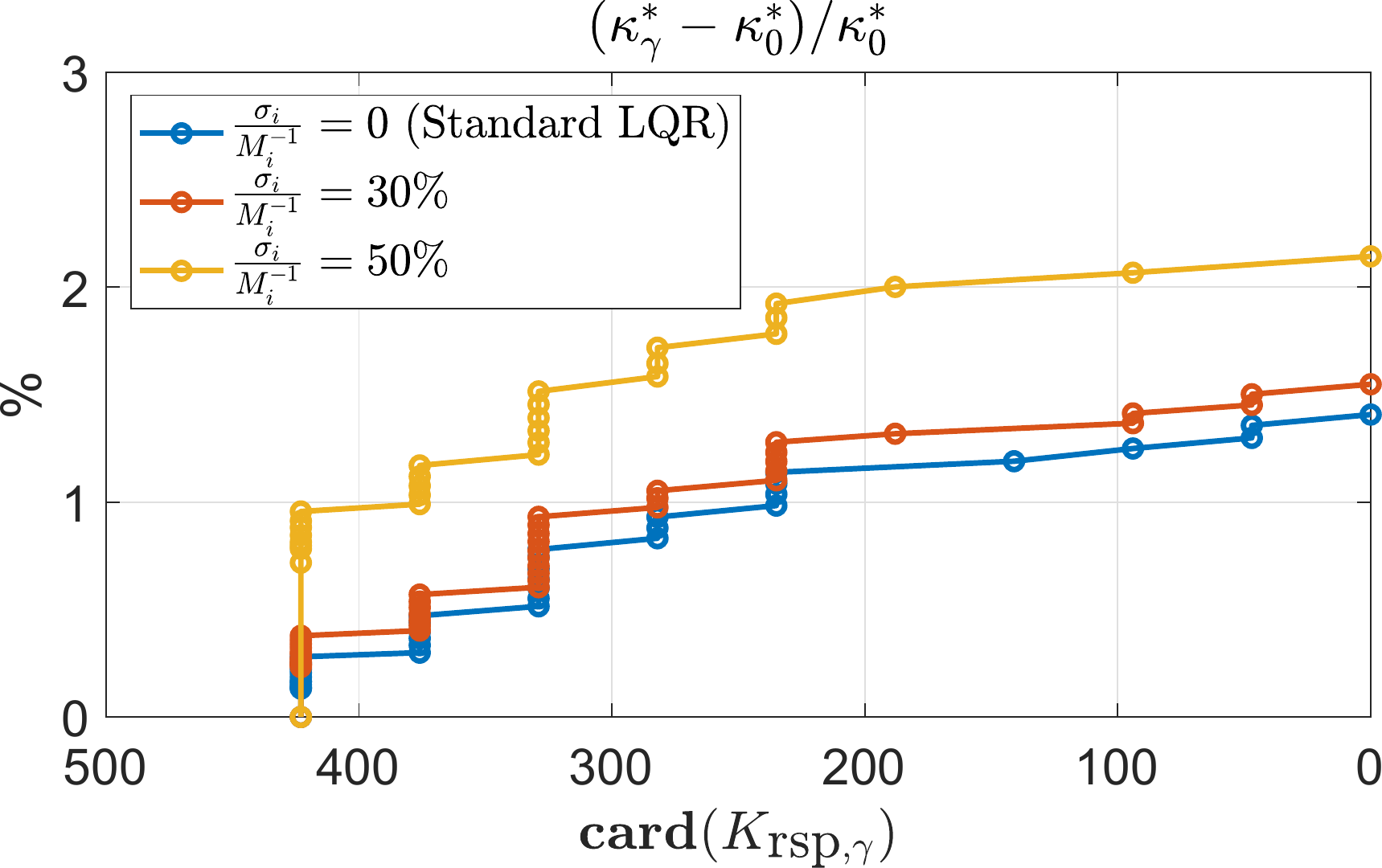}
\caption{Sparsity comparison under different row-sparsity regularizers. Note that increasing $\gamma$ results in more sparse pattern on $K_{\textrm{rsp},\gamma}$}
\label{fig:tradeoff_LQRm_Cost_kappa}
\end{figure}

The sparsity patterns for $\gamma = 0, 5~\textrm{and}~ 7$ are presented in Fig. \ref{fig:different_sparsity_K}. For $\gamma = 0$, the optimal feedback gain $K_0$ is fully populated, thereby requiring all 10 generators contributing to a mean-square stabilizing solution. As $\gamma$ increases, the rows of the state-feedback matrix $K_\gamma$ becomes sparse whereas the relative cost objective $(\kappa_\gamma^* - \kappa_0^*)/\kappa_0^*$ increases only slightly, see Fig. \ref{fig:tradeoff_LQRm_Cost_kappa}. In particular, for $\gamma = 5$, the identified control architecture indicates that the controllers of generators 2, 9 and 10 are not necessary to stabilize the system. As $\gamma$ increases to 7, most of stabilization burden is on generator 8 but only with 1.4\% LQRm cost increase.

\begin{figure}[htbp!]
\centering
\includegraphics[width=3.5in]{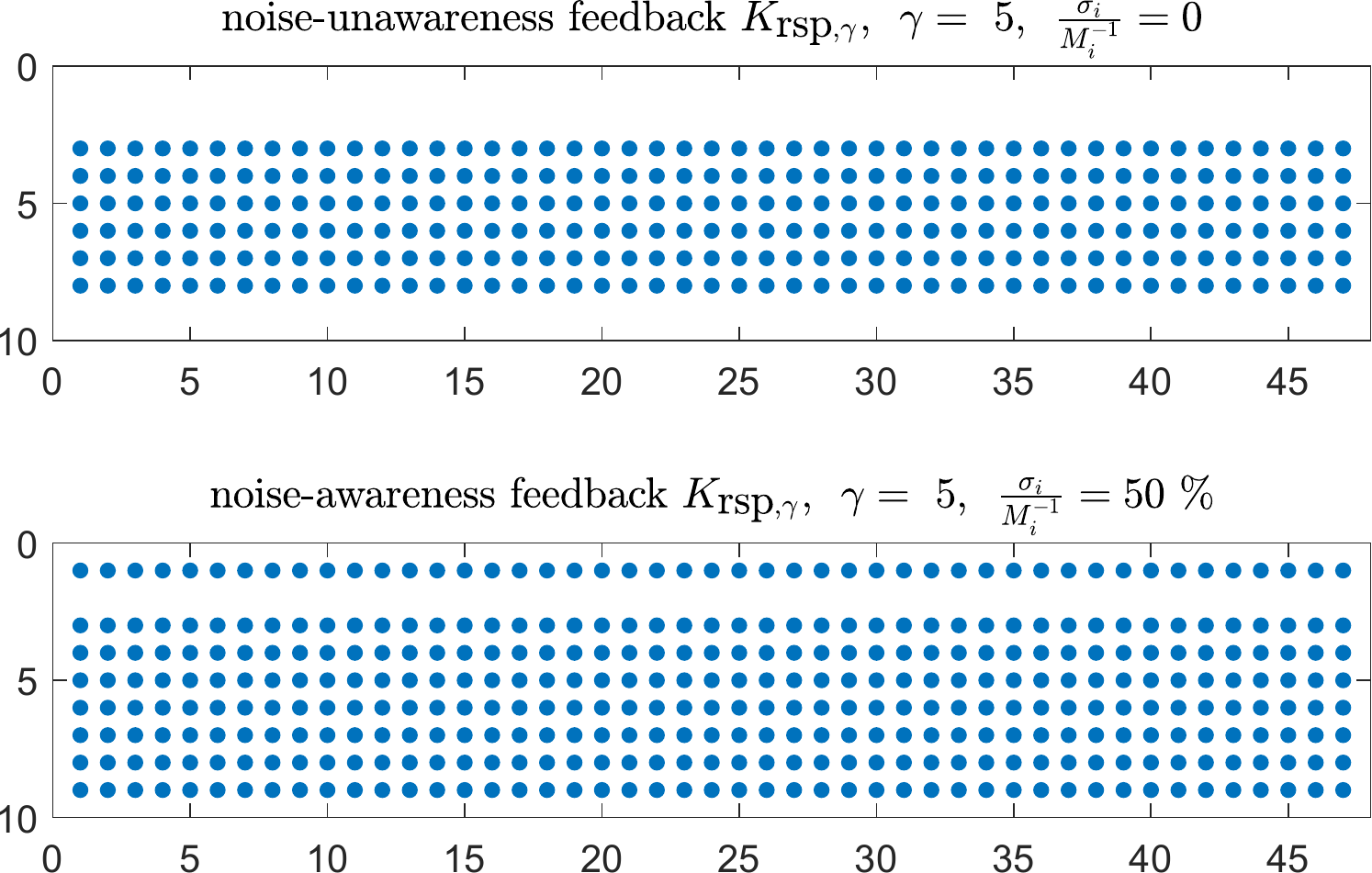}
\caption{Sparsity patterns of the noise-unaware and noise-aware state-feedback control law $K_{\textrm{rsp},\gamma}$.}
\label{fig:noise_awareness_Comparison}
\end{figure}

\begin{figure}[htbp!]
\centering
\includegraphics[width=3.5in]{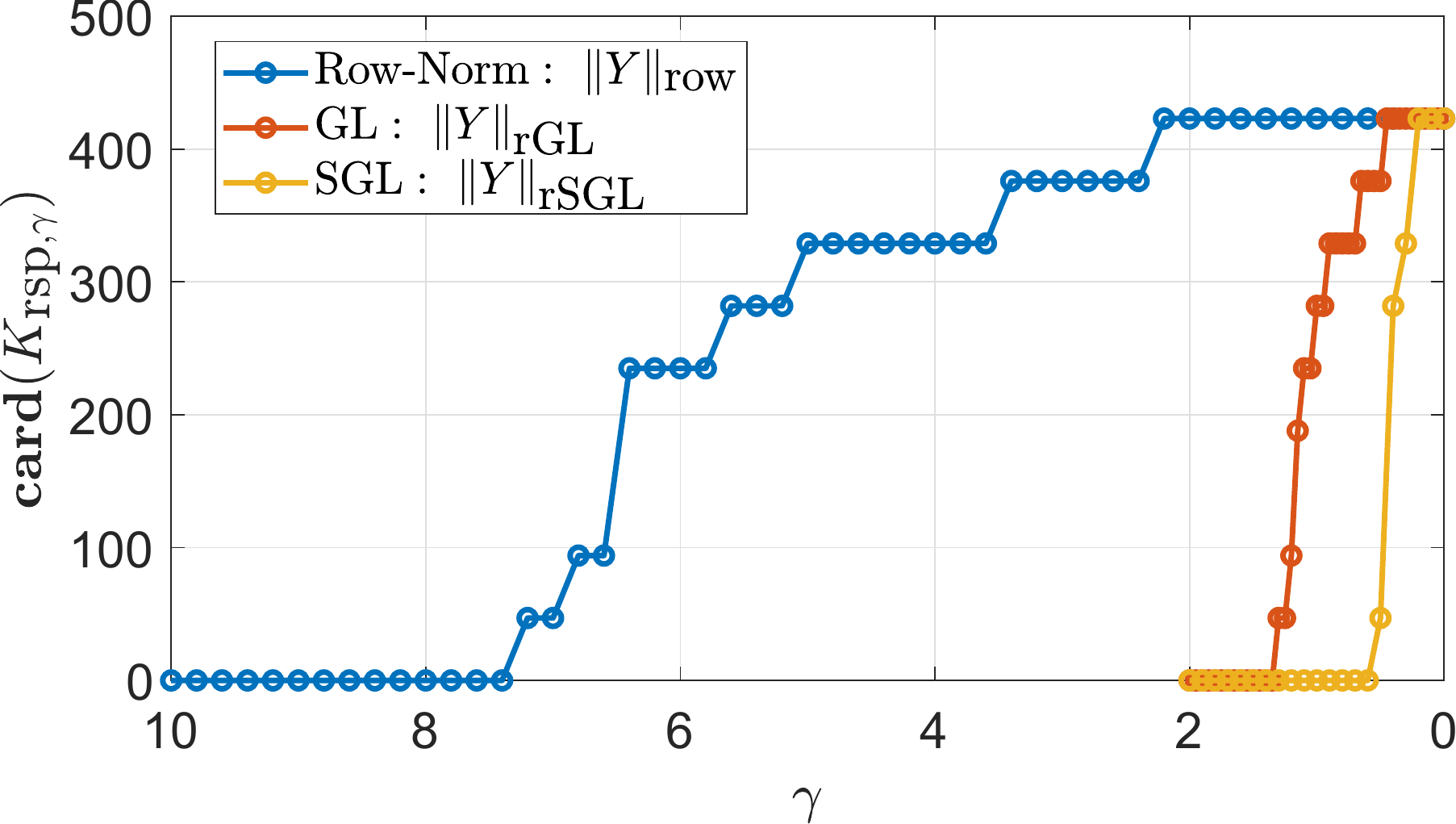}
\caption{Sparsity comparison under different row-sparsity regularizers with $\sigma_i/\overline{M}^{-1}_i= 10\%$.}
\label{fig:tradeoff_LQRm_Cost_gamma_different_regularizers}
\end{figure}

We next compare the LQRm cost with various levels of multiplicative noise as shown in Fig. \ref{fig:tradeoff_LQRm_Cost_kappa}. As the standard deviation $\sigma$ of the multiplicative noise increases, the LQRm cost increases since more control effort is required for stabilizing the system-level disturbances. In addition, the sparsity-promoting structures of the noise-aware ($\sigma_i/\overline{M}_i^{-1} = 50\%$) and noise-unaware ($\sigma = 0$)  state-feedback controllers $K_{\textrm{rsp},\gamma}$ are given in Fig. \ref{fig:noise_awareness_Comparison}. We observe that more generators need to be included for stabilizing the system due to the significant system-level disturbances. This also implies that a noise-unaware state-feedback controller may fail to stabilize a stochastic linear system with substantial multiplicative noise in a mean-square sense. In practice, the multiplicative noise is an inherent part of linearized system models and noisy control channels. This emphasizes the importance and necessity of having a noise-aware sparse architecture design to improve robustness to system-level disturbances for a mean-square stabilizing solution. Finally, we evaluate the proposed approach with various row sparsity induced norms, such as the group LASSO and sparse group LASSO, (see Fig. \ref{fig:tradeoff_LQRm_Cost_gamma_different_regularizers}). All of three sparsity-promoting norms successfully provide sparsity patterns for different $\gamma$. The sparse group LASSO (with $\mu = 0.5$) and the group LASSO lead to more aggressive sparsity patterns than the row-norm regularizer.

Overall, we conclude that the proposed approach successfully provides a sparse control solution, which reduces the number of controllers (with low-dimensional outputs) only at the expense of a small decrease of system performance. 


\section{Conclusions}
This paper proposes a sparse feedback control architecture design for stochastic linear systems with multiplicative noise. We minimize the sparsity-promoting matrix norms subject to a mean-square stability LMI condition as an SDP problem to approximate the nonconvex combinatorial problem. For a large-scale dynamic system with system-level disturbances, the designed sparse stabilizing solution successfully reduces the number of controllers, limits the unnecessary output information exchanges, and only slightly increases the LQRm cost.

\bibliography{reference}

\begin{thebibliography}{10}

\bibitem{harris1998signal}
C.~M. Harris and D.~M. Wolpert, ``Signal-dependent noise determines motor
  planning,'' {\em Nature}, vol.~394, no.~6695, p.~780, 1998.

\bibitem{todorov2005stochastic}
E.~Todorov, ``Stochastic optimal control and estimation methods adapted to the
  noise characteristics of the sensorimotor system,'' {\em Neural Computation},
  vol.~17, no.~5, pp.~1084--1108, 2005.

\bibitem{todorov2003optimal}
E.~Todorov and W.~Li, ``Optimal control methods suitable for biomechanical
  systems,'' in {\em 25th Annual International Conference of the IEEE
  Engineering in Medicine and Biology Society}, vol.~2, pp.~1758--1761, 2003.

\bibitem{todorov2002optimal}
E.~Todorov and M.~I. Jordan, ``Optimal feedback control as a theory of motor
  coordination,'' {\em Nature Neuroscience}, vol.~5, no.~11, p.~1226, 2002.

\bibitem{norris1989selection}
G.~Norris and R.~Skelton, ``Selection of dynamic sensors and actuators in the
  control of linear systems,'' 1989.

\bibitem{fardad2011sparsity}
M.~Fardad, F.~Lin, and M.~R. Jovanovi{\'c}, ``Sparsity-promoting optimal
  control for a class of distributed systems,'' in {\em Annual American Control
  Conference}, pp.~2050--2055, 2011.

\bibitem{schuler2011design}
S.~Schuler, P.~Li, J.~Lam, and F.~Allg{\"o}wer, ``Design of structured dynamic
  output-feedback controllers for interconnected systems,'' {\em International
  Journal of Control}, vol.~84, no.~12, pp.~2081--2091, 2011.

\bibitem{matni2017communication}
N.~Matni, ``Communication delay co-design in $\mathcal{H}_2$-distributed
  control using atomic norm minimization,'' {\em IEEE Transactions on Control
  of Network Systems}, vol.~4, no.~2, pp.~267--278, 2017.

\bibitem{fardad2014design}
M.~Fardad and M.~R. Jovanovic, ``On the design of optimal structured and sparse
  feedback gains via sequential convex programming,'' in {\em American Control
  Conference}, pp.~2426--2431, 2014.

\bibitem{matni2016regularization}
N.~Matni and V.~Chandrasekaran, ``Regularization for design,'' {\em IEEE
  Transactions on Automatic Control}, vol.~61, no.~12, pp.~3991--4006, 2016.

\bibitem{taha2018time}
A.~F. Taha, N.~Gatsis, T.~Summers, and S.~A. Nugroho, ``Time-varying sensor and
  actuator selection for uncertain cyber-physical systems,'' {\em IEEE
  Transactions on Control of Network Systems}, vol.~6, no.~2, pp.~750--762,
  2018.

\bibitem{summers2016actuator}
T.~Summers, ``Actuator placement in networks using optimal control performance
  metrics,'' in {\em 55th IEEE Conference on Decision and Control},
  pp.~2703--2708, 2016.

\bibitem{jovanovic2016controller}
M.~R. Jovanovi{\'c} and N.~K. Dhingra, ``Controller architectures: Tradeoffs
  between performance and structure,'' {\em European Journal of Control},
  vol.~30, pp.~76--91, 2016.

\bibitem{zare2018proximal}
A.~Zare, H.~Mohammadi, N.~K. Dhingra, T.~T. Georgiou, and M.~R. Jovanovi{\'c},
  ``Proximal algorithms for large-scale statistical modeling and
  sensor/actuator selection,'' {\em IEEE Transactions on Automatic Control},
  vol.~65, no.~8, pp.~3441--3456, 2019.

\bibitem{gravell2019sparse}
B.~Gravell, Y.~Guo, and T.~Summers, ``Sparse optimal control of networks with
  multiplicative noise via policy gradient,'' {\em IFAC-PapersOnLine}, vol.~52,
  no.~20, pp.~327--332, 2019.

\bibitem{summers2015submodularity}
T.~H. Summers, F.~L. Cortesi, and J.~Lygeros, ``On submodularity and
  controllability in complex dynamical networks,'' {\em IEEE Transactions on
  Control of Network Systems}, vol.~3, no.~1, pp.~91--101, 2015.

\bibitem{chang2018co}
C.-Y. Chang, S.~Mart{\'\i}nez, and J.~Cort{\'e}s, ``Co-optimization of control
  and actuator selection for cyber-physical systems,'' {\em IFAC-PapersOnLine},
  vol.~51, no.~23, pp.~118--123, 2018.

\bibitem{kim1991measure}
Y.~Kim and J.~L. Junkins, ``Measure of controllability for actuator
  placement,'' {\em Journal of Guidance, Control, and Dynamics}, vol.~14,
  no.~5, pp.~895--902, 1991.

\bibitem{dhingra2014admm}
N.~K. Dhingra, M.~R. Jovanovi{\'c}, and Z.-Q. Luo, ``An {ADMM} algorithm for
  optimal sensor and actuator selection,'' in {\em IEEE Annual Conference on
  Decision and Control}, pp.~4039--4044, 2014.

\bibitem{polyak2013lmi}
B.~Polyak, M.~Khlebnikov, and P.~Shcherbakov, ``An {LMI} approach to structured
  sparse feedback design in linear control systems,'' in {\em European Control
  Conference}, pp.~833--838, IEEE, 2013.

\bibitem{dhingra2016method}
N.~K. Dhingra and M.~R. Jovanovi{\'c}, ``A method of multipliers algorithm for
  sparsity-promoting optimal control,'' in {\em Annual American Control
  Conference}, pp.~1942--1947, 2016.

\bibitem{moghaddam2016customized}
S.~H. Moghaddam and M.~R. Jovanovi{\'c}, ``Customized algorithms for growing
  connected resistive networks,'' {\em IFAC-PapersOnLine}, vol.~49, no.~18,
  pp.~968--973, 2016.

\bibitem{boyd1994linear}
S.~Boyd, L.~El~Ghaoui, E.~Feron, and V.~Balakrishnan, {\em Linear Matrix
  Inequalities in System and Control Theory}.
\newblock SIAM, 1994.

\bibitem{wonham1967optimal}
W.~M. Wonham, ``Optimal stationary control of a linear system with
  state-dependent noise,'' {\em SIAM Journal on Control}, vol.~5, no.~3,
  pp.~486--500, 1967.

\bibitem{kleinman1969optimal}
D.~Kleinman, ``Optimal stationary control of linear systems with
  control-dependent noise,'' {\em IEEE Transactions on Automatic Control},
  vol.~14, no.~6, pp.~673--677, 1969.

\bibitem{gershon2001h}
E.~Gershon, U.~Shaked, and I.~Yaesh, ``$\mathcal{H}_{\infty}$ control and
  filtering of discrete-time stochastic systems with multiplicative noise,''
  {\em Automatica}, vol.~37, no.~3, pp.~409--417, 2001.

\bibitem{rami2001solvability}
M.~A. Rami, X.~Chen, J.~B. Moore, and X.~Y. Zhou, ``Solvability and asymptotic
  behavior of generalized {R}iccati equations arising in indefinite stochastic
  {LQ} controls,'' {\em IEEE Transactions on Automatic Control}, vol.~46,
  no.~3, pp.~428--440, 2001.

\bibitem{phillis1985controller}
Y.~Phillis, ``Controller design of systems with multiplicative noise,'' {\em
  IEEE Transactions on Automatic Control}, vol.~30, no.~10, pp.~1017--1019,
  1985.

\bibitem{arnold1974stochastic}
L.~Arnold, ``Stochastic differential equations,'' {\em New York}, 1974.

\bibitem{el1995state}
L.~El~Ghaoui, ``State-feedback control of systems with multiplicative noise via
  linear matrix inequalities,'' {\em Systems \& Control Letters}, vol.~24,
  no.~3, pp.~223--228, 1995.

\bibitem{ahsen2017two}
M.~E. Ahsen, N.~Challapalli, and M.~Vidyasagar, ``Two new approaches to
  compressed sensingexhibiting both robust sparse recovery and the grouping
  effect,'' {\em The Journal of Machine Learning Research}, vol.~18, no.~1,
  pp.~1745--1768, 2017.

\bibitem{rami2000linear}
M.~A. Rami and X.~Y. Zhou, ``Linear matrix inequalities, riccati equations, and
  indefinite stochastic linear quadratic controls,'' {\em IEEE Transactions on
  Automatic Control}, vol.~45, no.~6, pp.~1131--1143, 2000.

\bibitem{guo2019performance}
Y.~Guo and T.~H. Summers, ``A performance and stability analysis of low-inertia
  power grids with stochastic system inertia,'' in {\em Annual American Control
  Conference}, pp.~1965--1970, 2019.

\bibitem{hill1990stability}
D.~J. Hill and I.~M. Mareels, ``Stability theory for differential/algebraic
  systems with application to power systems,'' {\em IEEE Transactions on
  Circuits and Systems}, vol.~37, no.~11, pp.~1416--1423, 1990.

\bibitem{machowski2020power}
J.~Machowski, Z.~Lubosny, J.~W. Bialek, and J.~R. Bumby, {\em Power System
  Dynamics: Stability and Control}.
\newblock John Wiley \& Sons, 2020.

\bibitem{zimmerman2010matpower}
R.~D. Zimmerman, C.~E. Murillo-S{\'a}nchez, and R.~J. Thomas, ``{MATPOWER}:
  Steady-state operations, planning, and analysis tools for power systems
  research and education,'' {\em IEEE Transactions on Power Systems}, vol.~26,
  no.~1, pp.~12--19, 2010.

\bibitem{aps2019mosek}
M.~ApS, ``{MOSEK} optimization toolbox for {MATLAB},'' {\em User’s Guide and
  Reference Manual, Version}, vol.~4, 2019.

\bibitem{grant2014cvx}
M.~Grant and S.~Boyd, ``{CVX}: {MATLAB} software for disciplined convex
  programming, version 2.1,'' 2014.

\end{thebibliography}
\bibliographystyle{ieeetr}

\section*{Appendix: Results for Discrete-time Stochastic Linear Systems}\label{appendix}
In this appendix, we present the mean-square stability condition and the reformulation of the regularized LQRm problem for a discrete-time linear system with state- and input-dependent multiplicative noises. Consider a discrete-time stochastic linear system \cite{boyd1994linear}:
\begin{equation}\label{eq:discrete_dynamic}
\begin{aligned}
    x_{t+1} = \overline{A}_0x_t + \overline{B}_0u_t + \sum_{i=1}^k\overline{\sigma}_i\overline{A}_i x_tw_{it}  + \sum_{j=1}^\ell\overline{\rho}_j  
    \overline{B}_j u_tp_{jt},
\end{aligned}
\end{equation}
where $x_t\in \mathbb{R}^n$ denotes the state vector, $u_t\in\mathbb{R}^m$ denotes the control input vector. We use $w_{it}(i = 1,\ldots,k)$ and $p_{jt}(j = 1,\ldots,\ell)$ to denote the independent, identically random variables with
\begin{equation*}
\begin{aligned}
    &~~~~\mathbf{E}[w_{it}] = 0,~ \mathbf{E}[p_{jt}] = 0,~
    \mathbf{E}[w_{it}^2] = \overline{\sigma}_i^2,~
    \mathbf{E}[p_{jt}^2] = \overline{\rho}_j^2,\\
    &~~~~\mathbf{E}[w_{i_1,t}w_{i_2,t}] = 0~ (i_1 \neq i_2),~  \mathbf{E}[p_{j_1}p_{j_2}] = 0~ (j_1 \neq j_2),\\
    &\textrm{and} ~~~~~~~~~~~~~~~~~~~~~~~\mathbf{E}[w_{it}p_{jt}] = 0,~~~~~~~~~~~~~~~~~~~~~~~~~~~~\\
    &~~~~~~~~~~\forall i,i_1,i_2 = 1,\ldots,k, ~\textrm{and}~ \forall j,j_1,j_2 = 1,\ldots,\ell.~~~~~~~~
\end{aligned}
\end{equation*}
The scale factors $\overline{\sigma}_i$ and $\overline{\rho}_j$ indicate the standard deviation which normalize $w_{it}$ and $p_{jt}$ with the unit variance. The constant system matrices are $\overline{A}_0 \in \mathbb{R}^{n\times n}$ and $\overline{B}_0\in\mathbb{R}^{n \times m}$. The state dependent noise is allocated by matrix $\overline{A}_i\in \mathbb{R}^{n \times n}$, and the input dependent noise is allocated by matrix $\overline{B}_i \in \mathbb{R}^{n \times m}$. The system \eqref{eq:discrete_dynamic} is stabilizable via the state-feedback control $u_t = Kx_t$ if and only if there exists a matrix $X\in\mathbb{S}^n_{++}$ such that the following condition holds:
\begin{equation}\label{eq:discrete_stability_condition_1}\nonumber
\begin{aligned}
    (\overline{A}_0 + \overline{B}_0K)^\top X(\overline{A}_0 + \overline{B}_0K) - X + \sum_{i=1}^k \overline{\sigma}_i^2 \overline{A}_i^\top X\overline{A}_i \\
    + \sum_{j=1}^\ell \overline{\rho}_j^2 K^\top \overline{B}_j^\top X \overline{B}_j K \prec 0.
\end{aligned}
\end{equation}
We pre- and post-multiply the above inequality by $P = X^{-1}$ and introduce a new variable $Y = KP$, which leads to
\begin{equation}\label{eq:discrete_stability_condition_2}
\begin{aligned}
    (\overline{A}_0P + \overline{B}_0Y)^\top P^{-1}(\overline{A}_0P + \overline{B}_0Y) - P ~~~~~~~~~~~~~~~~~~~~\\
    + \sum_{i=1}^k \overline{\sigma}_i^2 P^\top\overline{A}_i^\top P^{-1}\overline{A}_iP + \sum_{j=1}^\ell \overline{\rho}_j^2 Y^\top \overline{B}_j^\top P^{-1} \overline{B}_j Y \prec 0.
\end{aligned}
\end{equation}
We then apply the Schur's Lemma on \eqref{eq:discrete_stability_condition_2} and come to a LMI:
\begin{equation} \label{eq:LMOstabilityDiscrete}
    \begin{bmatrix}
    P & (\overline{A}_0P + \overline{B}_0Y) & \overline{Z}\\
    (\overline{A}_0P + \overline{B}_0Y)^\top & P & 0\\
    \overline{Z}^\top & 0 & \overline{Z}_P
    \end{bmatrix} \succ 0,
\end{equation}
where $\overline{Z} = \begin{bmatrix} \overline{\sigma}_1 \overline{A}_1P, \ldots, \overline{\sigma}_k \overline{A}_kP, \overline{\rho}_1 \overline{B}_1Y,\ldots,\overline{\rho}_\ell \overline{B}_\ell Y \end{bmatrix}$ and 
\begin{equation*}
\overline{Z}_P = \textbf{blkdiag}\underbrace{\left(X, \ldots, X\right)}_{k+\ell}.
\end{equation*}
Equation \eqref{eq:LMOstabilityDiscrete} can replace the stability condition \eqref{eq:LMOstabilityContinuous} in the row sparsity-promoting problem \eqref{eq:rsp_problem} and the column sparsity-promoting problem \eqref{eq:csp_problem} when the system dynamic is given in the discrete-time domain. Similarly to \eqref{eq:LQR_SDP}, we present a relaxation of the regularized LQRm formulation for the stochastic discrete linear system \eqref{eq:discrete_dynamic} to minimize the upper bound of LQRm cost:
\begin{subequations}\label{eq:LQR_SDP_discrete}
    \begin{align}
    & \min_{Y,P, \kappa, \Pi} \quad \quad \kappa + \gamma \|Y\|_\textrm{reg},\\
    & \textrm{subject to:} \quad\quad \mathbf{Tr}(\Pi) \leq \kappa, \begin{bmatrix}
    \Pi & \Sigma^{\frac{1}{2}}_0\\
    \Sigma^{\frac{1}{2}}_0 & P
    \end{bmatrix} \succeq 0, P \succ 0, \\ 
    & \scalebox{0.95}[1]{$ \begin{bmatrix}
    P & (\overline{A}_0P + \overline{B}_0Y) & \overline{Z} & Y & P \\
    (\overline{A}_0P + \overline{B}_0Y)^\top & P & 0 & 0 & 0\\
    \overline{Z}^\top & 0 & \overline{Z}_P & 0 & 0\\
    Y^\top & 0 & 0 & R^{-1} & 0 \\
    P & 0 & 0 & 0 & Q^{-1} 
    \end{bmatrix} \succ 0, $}
    \end{align}
\end{subequations}
where $\kappa \in\mathbb{R}$ and $\Pi\in\mathbb{R}^{n\times n}$ are the slack variables. The importance of sparsity is defined by $\gamma$. The initial state condition is $\mathbf{E}[x_0x_0^\top] = \Sigma_0$ and the cost matrices $Q$ and $R$ are positive definite. The solution ${Y}^*$, ${P}^*$, ${\Pi}^*$ and ${\kappa}^*$ defines a sparse stabilizing law ${K}_{\textrm{sp}} = {Y}^*{P^*}^{-1}$ and the upper bound of the LQRm cost $J^*(K) = \mathbf{Tr}(\Sigma_0{P^*}^{-1})\leq {\kappa}^*$.

\end{document}